\documentclass[11pt,tbtags]{article}
\usepackage{amsmath,amssymb,amsthm,textcomp}
\usepackage[usenames]{color}\usepackage{graphicx}
\usepackage{enumerate}
\usepackage{url}

\numberwithin{equation}{section}

\allowdisplaybreaks
\newcounter{rot}

\definecolor{brown}{cmyk}{0, 0.72, 1, 0.45}
\definecolor{grey}{gray}{0.5}

\def\z{\zeta}

\allowdisplaybreaks
\def\hC{\hat{C}}

\def\ex{\mathbb{E\/}}

\def\part{\partial}
\def\tag#1 {\eqno(#1)}

\newcommand{\brac}[1]{\left( #1 \right)}
\newcommand{\bfrac}[2]{\left(\frac{#1}{#2}\right)}

  \def\d{\delta} 
    
  \def\k{\kappa} 
\def\z{\zeta}    \def\l{\lambda}
   \def\p{\pi}

\newtheorem{theorem}{Theorem}

\newtheorem{lemma}{Lemma}[section]
\newtheorem{sublemma}[lemma]{Sublemma}
\theoremstyle{remark}
\newtheorem{remark}{Remark}

\newcommand{\proofstart}{{\bf Proof\hspace{2em}}}
\newcommand{\proofend}{\hspace*{\fill}\mbox{$\Box$}}
\newcommand{\set}[1]{\left\{#1\right\}}
\newcommand{\ignore}[1]{}

\newcommand{\beq}[1]{\begin{equation}\label{#1}}
\def\eeq{\end{equation}}

\def\2G{{\sc 2greedy}}

\def\E{\mathbb{E}}

\title{On the length of a random minimum spanning tree.}
\date{March 9, 2013
}

\author{Colin Cooper\thanks{Department of  Computer Science, King's College,
University of London, London WC2R 2LS, UK} 
\and Alan Frieze\thanks{Department of Mathematical Sciences,
Carnegie Mellon University, Pittsburgh PA15217, USA.
Research supported in part by NSF Grant CCF1013110} 
\and Nate Ince\thanks{Department of Mathematical Sciences,
Carnegie Mellon University, Pittsburgh PA15217, USA}
\and Svante Janson\thanks{Department of Mathematics, Uppsala University,
  SE-75310 Uppsala, Sweden. 
Research  supported in part by the Knut and Alice Wallenberg Foundation} 
\and Joel Spencer\thanks{Courant Institute, New York, NY 10012 (USA)} 
}

\newenvironment{romenumerate}[1][0pt]{
\addtolength{\leftmargini}{#1}\begin{enumerate}
 }{\end{enumerate}}

\newcommand{\refT}[1]{Theorem~\ref{#1}}

\newcommand{\refL}[1]{Lemma~\ref{#1}}
\newcommand{\refR}[1]{Remark~\ref{#1}}

\begingroup
  \divide\count255 by 60
  \multiply\count255 by -60
  \advance\count255 by \time
  \ifnum \count255 < 10 \xdef\klockan{\the\count1.0\the\count255}
  \else\xdef\klockan{\the\count1.\the\count255}\fi
\endgroup

\newcommand{\sumj}{\sum_{j\ge1}}
\newcommand{\sumk}{\sum_{k\ge1}}

\newcommand\xpar[1]{(#1)}
\newcommand\bigpar[1]{\bigl(#1\bigr)}
\newcommand\Bigpar[1]{\Bigl(#1\Bigr)}
\newcommand\biggpar[1]{\biggl(#1\biggr)}
\newcommand\lrpar[1]{\left(#1\right)}

\newcommand\xcpar[1]{\{#1\}}

\newcommand\bigabs[1]{\bigl|#1\bigr|}

\newcommand\lrabs[1]{\left|#1\right|}
\def\rompar(#1){\textup(#1\textup)}    
\newcommand\xfrac[2]{#1/#2}

\newcommand\parfrac[2]{\lrpar{\frac{#1}{#2}}}

\newcommand\Bigparfrac[2]{\Bigpar{\frac{#1}{#2}}}

\def\xexp(#1){e^{#1}}
\newcommand\ceil[1]{\lceil#1\rceil}

\newcommand\ntoo{\ensuremath{{n\to\infty}}}

\newcommand\ktoo{\ensuremath{{k\to\infty}}}

\newcommand\punkt{.\spacefactor=1000}    
    
\newcommand\ie{i.e\punkt}
\newcommand\eg{e.g\punkt}

\newcommand\cf{cf\punkt}

\newcounter{CC}
\newcommand{\CC}{\stepcounter{CC}\CCx} 
\newcommand{\CCx}{C_{\arabic{CC}}}     
\newcommand{\CCdef}[1]{\xdef#1{\CCx}}     
\newcounter{cc}

\renewcommand\E{\operatorname{\mathbb E{}}}
\renewcommand\P{\operatorname{\mathbb P{}}}

\newcommand\Po{\operatorname{Po}}

\newcommand\gd{\delta}

\newcommand\gf{\varphi}

\newcommand\gG{\Gamma}

\newcommand\gl{\lambda}

\renewcommand\phi{\xxx}  

\newcommand\cB{\mathcal B}

\newcommand\ett[1]{\boldsymbol1\xcpar{#1}}

\newcommand\qw{^{-1}}

\newcommand\qqw{^{-1/2}}
\newcommand\qqq{^{1/3}}
\newcommand\qqqb{^{2/3}}

\newcommand\qqqbw{^{-2/3}}

\renewcommand{\=}{:=}

\newcommand\oi{[0,1]}

\newcommand\oooo{(-\infty,\infty)}

\newcommand\gnp{\ensuremath{G_{n,p}}}

\newcommand\kkc{\kk_{\mathsf c}}
\newcommand\kkcn{\kk_{\mathsf c, n}}
\newcommand\kk{\kappa}
\newcommand\pn{n^{-1}+\gl n^{-4/3}}
\newcommand\wl{w_\ell}
\newcommand\wj{w_{j+1}}

\newcommand\intp{\int_{p=0}^1}
\newcommand\ints{\int_{s=0}^1}

\newcommand\intx{\int_{x=0}^\infty}
\newcommand\inty{\int_{y=0}^\infty}
\newcommand\intgl{\int_{\gl=-n^{1/3}}^{n^{4/3}-n^{1/3}}}
\newcommand\intgloo{\int_{\gl=-\infty}^{\infty}}
\newcommand\ff{\bar f}
\newcommand\kx{\ceil{xn\qqqb}}
\newcommand\ccc{c_{2c}}
\newcommand\cca{c_{2a}}
\newcommand\ccb{c_{2b}}
\newcommand\bex{\cB_{\mathrm{ex}}}
\newcommand\fbex{f_{\mathrm{ex}}}

\begin{document}

\maketitle
\begin{abstract}
We study the expected value of the length $L_n$ of the minimum spanning tree of the complete graph $K_n$ when
each edge $e$ is given an independent uniform $[0,1]$ edge weight. We sharpen the result of Frieze
\cite{F1} that $\lim_{n\to\infty}\E(L_n)=\z(3)$ and show that 
$\E(L_n)=\z(3)+\frac{c_1}{n}+\frac{c_2+o(1)}{n^{4/3}}$ where
$c_1,c_2$ are explicitly defined constants.
\end{abstract}
\section{Introduction}
We study the expected value of the length $L_n$ of the minimum spanning tree of the complete graph $K_n$ when
each edge $e$ is given an independent uniform $[0,1]$ edge weight $X_e$. It was shown in Frieze \cite{F1} that
\beq{E1}
\lim_{n\to\infty}\E(L_n)=\z(3)=\sum_{k=1}^\infty \frac{1}{k^3}=1.202\ldots
\eeq
Since then there have been several generalisations and improvements. Steele \cite{Ste1} extended the applicability of
\eqref{E1} distribution-wise. Janson \cite{Ja} proved a central limit theorem for $L_n$. Penrose \cite{Pe},
Frieze and McDiarmid \cite{FM},
Beveridge, Frieze and McDiarmid \cite{BFM}, Frieze, Ruszink\'o and Thoma
\cite{FRT} analysed $L_n$ for graphs 
other than the complete graph.
Fill and Steele \cite{Ste2} used the Tutte polynomial to compute $\E(L_n)$
exactly for small values and 
Gamarnik \cite{Ga} computed $\E_{exp}(L_n)$ exactly
up to $n\leq 45$ using a more efficient algorithm,
where $\E_{exp}(L_n)$ is the expectation when 
the distribution of the $X_e$ is exponential with mean one. 
Li and Zhang \cite{LZ} consider more general distributions and prove in 
particular that
\beq{eqLZ}
\E_{exp}(L_n)-\E(L_n)=\frac{\z(3)}{n}+O\bfrac{\log^2n}{n^2}.
\eeq
Flaxman \cite{Flax} gives an upper bound on the lower tail of $L_n$.

Equation \eqref{E1} says that $\E(L_n)=\z(3)+o(1)$ as $n\to\infty$. Ideally, one would like to have an exact expansion
for $\E(L_n)$ as there is for the assignment problem, see W\"astlund
\cite{Wastlund} and the references therein. 
Such an expansion has proven elusive.  In this work we
improve the asymptotics of $E[L_n]$ by giving the secondary
and tertiary terms.
\begin{theorem}\label{th1}
$$\E(L_n)=\z(3)+\frac{c_1}{n}+\frac{c_2+o(1)}{n^{4/3}}$$
where
\begin{align*}
c_1&=-1-\zeta(3)-
\frac12\int_{x=0}^\infty \log\bigl(1-(1+x)e^{-x}\bigr) \, dx
\intertext{and}
c_2&
=\intx
\Bigpar{x^{-3}\psi(x^{3/2})e^{-x^3/24}-x^{-3}-\sqrt{\frac{\pi}8}x^{-3/2}
-\frac12}\,dx
\\
&=\frac23\inty
\Bigpar{y^{-2}\psi(y)e^{-y^2/24}-y^{-2}-\sqrt{\frac{\pi}8}y^{-1}
-\frac12}y^{-1/3}\,dy
\end{align*}
with $\psi$ defined in \eqref{mgf} below.
\end{theorem}
The two integral expressions defining $c_2$ are equal by the change of
variable $x=y^{2/3}$. 

A numerical integration (with \texttt{Maple}) yields
$c_1=0.0384956\dots$.
This shows that the rate of convergence to $\z(3)$ is order $1/n$ and is
from above. 
Further numerical computations show that 
$  c_2\approx -1.7295$,   
and these are explained
in an appendix.

To define $\psi$, we let the random variable 
$\bex=\ints B_{\mathrm{ex}}(s)\,ds$ 
be the area under a normalized Brownian excursion; we then let
\begin{equation}
  \label{mgf}
\psi(t)=\E e^{t \bex},
\end{equation}
the moment generating function $\psi$ of $\bex$.
The Brownian excursion area $\bex$ and its moments $\E \bex^\ell$ and moment
generating function $\psi$ have
been studied by several authors, see \eg{} Louchard \cite{Lou0,Lou}
and the survey by Janson \cite{SJ201}, where further references are given.
From these results, we derive an expression, see \eqref{jk}, that will show
$c_2$  
is well-defined. 
Note that $\psi(t)$ is finite for all $t>0$ (and thus \eqref{mgf}
holds for all complex $t$); 
indeed, 
see \cite[(53)]{SJ201} and the references there,
it is well-known that
\begin{equation}\label{ebex}
\E\bex^\ell\sim\sqrt{18}\,\ell\,(12e)^{-\ell/2}\ell^{\ell/2}  
\qquad \text{as $\ell\to\infty$},
\end{equation}
and thus \cite[Lemma 4.1(ii)]{SJISE} implies, 
\cf{} \cite[Remarks 3.1 and 4.9]{SJISE} (where  $\xi=2\bex$),
\begin{equation}
\label{jesp}
\psi(t) 
\sim
\tfrac12{t^2}e^{\xfrac{t^2}{24}}
\qquad \text{as }
t\to+\infty.
\end{equation}
More precisely, Janson and Louchard \cite{SJ203} show  
that the density $\fbex$ of $\bex$ satisfies
\begin{equation}
  \label{fbex}
\fbex(x)=\frac{72\sqrt6}{\sqrt\pi}x^2e^{-6x^2}\bigpar{1+O(x^{-2})},
\qquad x>0,
\end{equation}
from which routine calculations show that
\begin{equation}\label{jk}
  \psi(t)=\intx e^{tx}\fbex(x)\,dx
=\frac{t^2}2e^{t^2/24} \bigpar{1+O(t^{-2})},
\qquad t>0.
\end{equation}
Hence the integrand in the second integral defining $c_2$ in \refT{th1} 
is $O(y^{-4/3})$ as $y\to\infty$. Moreover, $\psi(0)=1$ and
$\psi'(0)=\E\bex=\sqrt{\pi/8}$, and thus a Taylor expansion shows that the
integrand is $O(y^{-1/3})$ as $y\to0$. 
(Similarly, the integrand in the first integral is $O(x^{-3/2})$ and $O(1)$.)
Consequently, the integrals defining $c_2$ converge absolutely.

\section{Proof of Theorem \ref{th1}}
We prove the theorem by using the expression (see Janson \cite{Ja}),
\beq{E2}
\E(L_n)=\int_{p=0}^1\E(\k(G_{n,p}))dp -1.
\eeq
Here $\k(G_{n,p})$ is the (random) number
of components in the random graph $G_{n,p}$. 

To evaluate \eqref{E2} we let $\k(k,j,p)=\k_n(k,j,p)$
denote the number of components of $G_{n,p}$ with $k$ vertices and $k+j$
edges in $G_{n,p}$.
The components neatly split into three categories:  trees ($j=-1$),
unicyclic ($j=0$) and complex ($j\geq 1$) components.  These
are evaluated separately.
\begin{lemma}\label{lem1}\
\begin{enumerate}[(a)]
 \item 
$$
\int_{p=0}^1\sum_{k\geq 1}\E(\k(k,-1,p))dp=\z(3)+\frac{3(\z(2)-\z(3))}{2n}
-\frac{1}{n^{4/3}}\int_{x=0}^\infty x^{-3}(1-e^{-x^3/24})\,dx+o(n^{-4/3}).
$$
\item 
\begin{multline*}
\int_{p=0}^1\sum_{k\geq3}\E(\k(k,0,p))dp=\frac{1}{2n}\brac{\z(3)-3\z(2)-
\int_{x=0}^\infty \log\bigl(1-(1+x)e^{-x}\bigr) \, dx.}\\
-\frac{\sqrt{\p/8}}{n^{4/3}}\int_{x=0}^\infty x^{-3/2}(1-e^{-x^3/24})\,dx
+o(n^{-4/3}).
\end{multline*}
\item 
With $\psi_2(x)=\psi(x)-1-\sqrt{\pi/2}\,x$,
	\begin{equation*}
\int_{p=0}^1\sum_{k\geq 1}\sum_{j\geq 1}\E(\k(k,j,p))dp=
1-\frac1n+\frac1{n^{4/3}}\intx \Bigpar{x^{-3}\psi_2(x^{3/2})e^{-x^3/24}-\frac12}\,dx
+o(n^{-4/3}).	  
	\end{equation*}
\end{enumerate}
\end{lemma}
\begin{remark}
Tree components contribute the main $\zeta(3)$ addend.
Unicyclic components contribute a secondary $O(\frac{1}{n})$
addend.  Roughly speaking there are no complex components for
$p\leq \frac{1}{n}$ and precisely one complex component (the
famous ``giant component") for $p\geq \frac{1}{n}$.  Were this
to be precisely the case the contribution of complex components
would be $1-\frac{1}{n}$.  The additional $\Theta(n^{-4/3})$
term in Lemma \ref{lem1} (c) comes from the behavior of complex
components in the critical window $p = \frac{1}{n} + \lambda n^{-4/3}$.
\end{remark}

\begin{remark} \label{Rc2ab}
The coefficients of $n^{-4/3}$ in \refL{lem1}(a) and (b) are easily evaluated
  as
$-\frac18 3^{-2/3} \Gamma(1/3)$  
and 
$-\frac123^{-1/6}\sqrt\pi\,\Gamma(5/6)$,   
respectively, see the appendix.  
The coefficient in (c) is expressed as an infinite sum and evaluated
numerically in the appendix. 
\end{remark}

\proofstart
We assume in the proof tacitly that $n$ is large enough when necessary.
We let $C_1,\dots$ denote some unimportant universal constants.

Let
$C(k,\ell)$ be the number of connected graphs on vertex set $[k]$ with
$\ell$ edges.
We begin by noting the standard formula
\begin{equation}
  \label{ekk}
\E \k(k,j,p) 
=\binom{n}{k}C(k,k+j)p^{k+j}(1-p)^{k(n-k)+\binom{k}{2}-k-j}.
\end{equation}
By Cayley's formula,
$C(k,k-1)=k^{k-2}$. Moreover, Wright \cite{Wright} proved that for every fixed
$j\ge-1$, 
\begin{equation}\label{wright}
  C(k,k+j)\sim \wj k^{k+3j/2-1/2}
\qquad\text{as \ktoo},
\end{equation}
for some constants $\wl>0$.
(See also \cite[\S8]{SJbirth} and the references there.
In the notation of \cite{Wright}, $w_{j+1}=\rho_j$.) 
We have $w_0=1$ and $w_1=\sqrt{\pi/8}$.
It was shown in Spencer \cite{Spencer} that 
\begin{equation}\label{wbex}
  \wl=\frac{\E \bex^\ell}{\ell!},
\qquad \ell\ge0,
\end{equation}
where $\bex$ is the Brownian excursion area defined above.
See further Janson \cite{SJ201}. 
Hence,
\begin{equation}
  \label{mgf2}
\psi(t)=\E e^{t \bex}=\sum_{\ell=0}^\infty \wl t^\ell.
\end{equation}

 Let
\begin{equation}\label{A}
  \begin{split}
A(k,k+j)&=\int_{p=0}^1\ex(\k(k,j,p))\,dp
\\
&=\binom{n}{k}C(k,k+j)\int_{p=0}^1p^{k+j}(1-p)^{k(n-k)+\binom{k}{2}-k-j}\,dp
\\
&=\binom{n}{k}C(k,k+j)\frac{(k+j)!\,(k(n-k)+\binom{k}{2}-k-j)!}{(k(n-k)+\binom{k}{2}+1)!}
\\
&=\frac{C(k,k+j)\,(k+j)!}{k!} \times B(k,k+j)
  \end{split}
\end{equation}
where, provided $k\le n$ and $k+j\le \binom k2$ (as in our case),
\begin{equation}\label{B}
  \begin{split}
B(k,k+j)&=\frac{n!}{(n-k)!}\cdot
\frac{(k(n-k)+\binom{k}{2}-k-j)!}{(k(n-k)+\binom{k}{2}+1)!}
\\
&=\frac{1}{n^{j+1}k^{k+j+1}}\frac{\prod_{i=0}^{k-1}\brac{1-\frac{i}{n}}}{\prod_{i=0}^{k+j}
\brac{1-\frac{k+1}{2n}-\frac{i-1}{kn}}}
\\
&=\frac{1}{n^{j+1}k^{k+j+1}}\exp\set{\sum_{m=1}^\infty\frac{1}{mn^m}\brac{\sum_{i=0}^{k+j}
\brac{\frac{k+1}{2}+\frac{i-1}{k}}^m-\sum_{i=0}^{k-1}i^m}}
\\
&=\frac{1}{n^{j+1}k^{k+j+1}}\exp\set{\sum_{m=1}^\infty\frac{t_{m}(k,j)}{mn^m}}.	
  \end{split}
\end{equation}
Observe that  as $\sum_{i=1}^a i^m \ge \int_0^a x^m \,dx$, for $\ell=k+j$  we have
\begin{align}\label{0}
t_m(k,j) &=\sum_{i=0}^\ell
\brac{\frac{k+1}{2}+\frac{i-1}{k}}^m
-\sum_{i=0}^{k-1}i^m
\leq(\ell+1)\brac{\frac{k+1}{2}+\frac{\ell-1}{k}}^m
-\frac{(k-1)^{m+1}}{m+1}.
\end{align}
This implies that, as is easily verified, 
\beq{00}
t_m(k,j)\leq 0\text{ if }m\geq2\text{ and }j\in\{0,-1\}\text{ and }k\geq 100.
\eeq

{\bf Case (a):} $1\leq k\leq n$, $j=-1$ (Tree components).

Now  we have by \eqref{B}
\begin{align*}
B(k,k-1)&=\frac{1}{k^{k}}\exp\left\{\frac{1}{n}\sum_{i=0}^{k-1}\brac{\frac{k+1}{2}+\frac{i-1}{k}}-
\frac{1}{n}\sum_{i=0}^{k-1}i\right.
\\
&\hspace{.84in}+\frac{1}{2n^2}\sum_{i=0}^{k-1}\brac{\frac{k+1}{2}+\frac{i-1}{k}}^2-\frac{1}{2n^2}\sum_{i=0}^{k-1}i^2
+\xi\Bigr\}
\end{align*}
where, using \eqref{00}, 
\begin{align}
&|\xi|\leq \sum_{m=3}^\infty\frac{10^{2m+1}}{mn^m}=O(n^{-3})&1\leq k\leq 100,
\label{100}\\
&0\geq \xi\geq -\sum_{m=3}^\infty\frac{k^{m+1}}{m(m+1)n^m}\geq
-\frac{k^4}{n^3}&k>100,\label{1001}
\end{align}
and hence for all $k\le n$,
\begin{equation}\label{1002} 
\xi=O(k^4/n^3). 
\end{equation}

This implies, after some calculation, that, for $1\le k\le n$,
$$B(k,k-1)=\frac{1}{k^{k}}\exp\set{\frac{3(k-1)}{2n}-\frac{k^3}{24n^2}+
O\brac{\frac{k^2}{n^2}+\frac{k^4}{n^3}}}$$
and then, by \eqref{A},
\begin{align*}
\sum_{k=1}^{n^{0.7}}A(k,k-1)&=\sum_{k=1}^{n^{0.7}}\frac{k^{k-2}}{k}\cdot  B(k,k-1)\\
&=\sum_{k=1}^{n^{0.7}}
\frac{1}{k^3} \exp\set{\frac{3(k-1)}{2n}-\frac{k^3}{24n^2}+O\brac{\frac{k^2}{n^2}+\frac{k^4}{n^3}}}\\
&=\sum_{k=1}^{n^{0.7}}\frac{e^{-k^3/24n^2}}{k^3}\brac{1+\frac{3(k-1)}{2n}+O\brac{\frac{k^2}{n^2}+\frac{k^4}{n^3}}}.
\end{align*}
Now, by simple estimates, 
\beq{x1a}
\sum_{k=1}^{n^{0.7}}\frac{e^{-k^3/24n^2}}{k^3}\times 
O\brac{\frac{k^2}{n^2}+\frac{k^4}{n^3}}=O(n^{-5/3})
\eeq
and
\begin{equation}\label{x2a}  
\begin{split}
\sum_{k=1}^{n^{0.7}}\frac{(1-e^{-k^3/24n^2})}{k^3}\brac{1+\frac{3(k-1)}{2n}}
&=o(n^{-4/3})+\sum_{k=n^{2/3}/\ln n}^{n^{2/3}\ln n}
\frac{(1-e^{-k^3/24n^2})}{k^3}\\
&=o(n^{-4/3})+\frac{1}{n^{4/3}}\int_{x=0}^\infty x^{-3}(1-e^{-x^3/24})\,dx.
\end{split}  
\end{equation}
Thus
\begin{align}
&\sum_{k=1}^{n^{0.7}} A(k,k-1)
\nonumber\\
&\quad=\sum_{k=1}^{n^{0.7}}\frac{1}{k^3}+\frac{1}{n}\sum_{k=1}^{n^{0.7}} \frac{3(k-1)}{2k^3}
-\frac{1}{n^{4/3}}\int_{x=0}^\infty x^{-3}(1-e^{-x^3/24})\,dx+o(n^{-4/3})
\nonumber\\
&\quad=\z(3)+O(n^{-1.4})+\frac{3(\z(2)-\z(3))}{2n}+O(n^{-1.7})
-\frac{1}{n^{4/3}}\int_{x=0}^\infty x^{-3}(1-e^{-x^3/24})\,dx+o(n^{-4/3})
\nonumber\\
&\quad=\z(3)+\frac{3(\z(2)-\z(3))}{2n}-
\frac{1}{n^{4/3}}\int_{x=0}^\infty x^{-3}(1-e^{-x^3/24})\,dx+o(n^{-4/3}).	
\label{A1}
\end{align}
When $k\geq n^{0.7}$ we have from \eqref{B} and \eqref{00} that 
\begin{equation*}
B(k,k-1)\le\frac{1}{k^{k}}\exp\lrpar{\frac{1}{n}
\sum_{i=0}^{k-1}\brac{\frac{k+1}{2}+\frac{i-1}{k}}
- \frac{1}{n}\sum_{i=0}^{k-1}i}
=
\frac{1}{k^{k}}\exp\set{\frac{3(k-1)}{2n}}
\le \frac{e^{3/2}}{k^k}.
\end{equation*}

This implies that $A(k,k-1)\leq k^{-3}e^{3/2}$.
This gives
$$\sum_{k> n^{0.7}}A(k,k-1)\leq \sum_{k> n^{0.7}}\frac{e^{3/2}}{k^3}
=O(n^{-1.4})=o(n^{-4/3}).$$
Together with \eqref{A1}, this verifies (a).

{\bf Case (b):} $1\leq k\leq n$, $j=0$ (Unicyclic components).

R\'enyi \cite{Re} proved 
(see \eg{} Bollob\'as \cite[Theorem 5.18]{B})
that, cf.\ the more general \eqref{wright} above,
\beq{renyi}
C(k,k)=\frac{(k-1)!}{2}\sum_{l=0}^{k-3}\frac{k^l}{l!}\sim \sqrt{\frac{\pi}{8}}k^{k-1/2}.
\eeq
Now for $1\le k\le n$ we have by \eqref{B}
\begin{align*}
B(k,k)=\frac{1}{nk^{k+1}}
\exp&\left\{\frac{1}{n}\sum_{i=0}^k\brac{\frac{k+1}{2}+\frac{i-1}{k}}-
\frac{1}{n}\sum_{i=0}^{k-1}i\right.\\
&\qquad
\left.+\frac{1}{2n^2}\sum_{i=0}^k\brac{\frac{k+1}{2}+\frac{i-1}{k}}^2-\frac{1}{2n^2}\sum_{i=0}^{k-1}i^2
+\xi\right\}
\end{align*}
where $\xi$ satisfies \eqref{100}--\eqref{1002}.
Thus, after some calculation,
$$B(k,k)=\frac{1}{k^{k+1}n}\exp\set{\frac{2k}{n}-\frac{1}{kn}-\frac{k^3}{24n^2}
+O\brac{\frac{k^2}{n^2}+\frac{k^4}{n^3}}}$$
and then
\begin{equation}\label{A44}
  \begin{split}
\sum_{k=3}^{n^{0.7}}A(k,k)&=\frac{1}{n}\sum_{k=3}^{n^{0.7}}\frac{C(k,k)}{k^{k+1}}
\exp\set{-\frac{k^3}{24n^2}+
O\brac{\frac{k}{n}+\frac{k^4}{n^3}}}\\
&=\frac{1}{n}\sum_{k=3}^{n^{0.7}}\frac{C(k,k)e^{-k^3/24n^2}}{k^{k+1}}
\set{1+O\brac{\frac{k}{n}+\frac{k^4}{n^3}}}.	
  \end{split}
\end{equation}
Now \eqref{renyi} implies
\beq{x1}
\frac{1}{n}\sum_{k=3}^{n^{0.7}}\frac{C(k,k)e^{-k^3/24n^2}}{k^{k+1}}\times 
O\brac{\frac{k}{n}+\frac{k^4}{n^3}}=O(n^{-5/3})
\eeq
and
\begin{equation}\label{x2}
  \begin{split}
\frac{1}{n}\sum_{k=3}^{n^{0.7}}\frac{C(k,k)(1-e^{-k^3/24n^2})}{k^{k+1}}
&=o(n^{-4/3})+\frac1n\sum_{k=n^{2/3}/\ln n}^{n^{2/3}\ln n}
\frac{C(k,k)(1-e^{-k^3/24n^2})}{k^{k+1}}\\
&=o(n^{-4/3})+\frac{\sqrt{\p/8}}{n^{4/3}}\int_{x=0}^\infty
x^{-3/2}(1-e^{-x^3/24})\,dx.
  \end{split}
\end{equation}
It follows from \eqref{A44}, \eqref{x1} and \eqref{x2} that
\beq{A5}
\sum_{k=3}^{n^{0.7}}A(k,k)=\frac{1}{n}\sum_{k=3}^\infty\frac{C(k,k)}{k^{k+1}}-
\frac{\sqrt{\p/8}}{n^{4/3}}\int_{x=0}^\infty x^{-3/2}(1-e^{-x^3/24})\,dx
+o(n^{-4/3})
.\eeq

For $k>n^{0.7}$ we observe that $t_1(k,0)\leq 2k$ in \eqref{0} and $t_m(k,0)\leq 0$ for $m\geq 2$ and so
$$B(k,k)\leq\frac{e^{2}}{k^{k+1}n}$$
and so
$$A(k,k)\leq e^{2}\frac{C(k,k)}{k^{k+1}n}=O\bfrac{1}{k^{3/2}n}.$$
It follows from this that
\beq{xx1}
\sum_{k=n^{0.7}}^nA(k,k)=O(n^{-1.35})=o(n^{-4/3}).
\eeq
We are almost done, we need to simplify the sum $\sum_{k=3}^\infty\frac{C(k,k)}{k^{k+1}}$.

Now, by \eqref{renyi},
\begin{equation}\label{c1}
\sum_{k=3}^\infty \frac{2C(k,k)}{k^{k+1}}
= \sum_{k=3}^\infty \frac{(k-1)!}{k^{k+1}}\sum_{i=0}^{k-3} \frac{k^i}{i!}
=
\sum_{i=0}^\infty \sum_{k=i+3}^\infty \frac{k^i}{k^{k+1}}\frac{(k-1)!}{i!}.
\end{equation}
In the last double sum, let us also add the terms with $k=i+2$, $k=i+1$ and
$k=i\ge1$. 
The terms with $k=i+2$ add up to
\begin{equation*}
  \begin{split}
 \sum_{k=2}^\infty \frac{k^{k-2}}{k^{k+1}}\frac{(k-1)!}{(k-2)!}
=
 \sum_{k=2}^\infty \frac{k-1}{k^{3}}
=
 \sum_{k=1}^\infty \frac{k-1}{k^{3}}
=\zeta(2)-\zeta(3).
  \end{split}
\end{equation*}
The terms with $k=i+1$ add up to
\begin{equation*}
  \begin{split}
 \sum_{k=1}^\infty \frac{k^{k-1}}{k^{k+1}}\frac{(k-1)!}{(k-1)!}
=
 \sum_{k=1}^\infty \frac{1}{k^{2}}
=\zeta(2).
  \end{split}
\end{equation*}
The terms with $k=i\ge1$ add up to
\begin{equation*}
  \begin{split}
 \sum_{k=1}^\infty \frac{k^{k}}{k^{k+1}}\frac{(k-1)!}{k!}
=
 \sum_{k=1}^\infty \frac{1}{k^{2}}
=\zeta(2).
  \end{split}
\end{equation*}
Consequently, \eqref{c1} yields
\beq{c1.2}
\sum_{k=3}^\infty \frac{2C(k,k)}{k^{k+1}}
=\zeta(3)-3\zeta(2)
+\sum_{k=1}^\infty\sum_{i=0}^k	\frac{k^i}{k^{k+1}}\frac{(k-1)!}{i!}\\
= \zeta(3)-3\zeta(2)
+\sum_{k=1}^\infty\sum_{i=0}^k	\frac{k!}{i!} k^{i-k-2}.
\eeq
We transform the sum further:
\begin{equation*}
  \begin{split}
\sum_{k=1}^\infty\sum_{i=0}^k	\frac{k!}{i!} k^{i-k-2}
&=
\sum_{k=1}^\infty\sum_{i=0}^k	\binom ki (k-i)!\, k^{i-k-2}
\\
&=
\sum_{k=1}^\infty\sum_{i=0}^k	\binom ki k^{-1} \int_{x=0}^\infty x^{k-i}
e^{-kx} \, dx
\\
&=
\int_{x=0}^\infty \sum_{k=1}^\infty\sum_{i=0}^k	k^{-1} \binom ki 
x^{k-i} e^{-kx} \, dx
\\
&=
\int_{x=0}^\infty \sum_{k=1}^\infty	k^{-1} (1+x)^{k} e^{-kx} \, dx
\\
&=
\int_{x=0}^\infty -\log\bigl(1-(1+x)e^{-x}\bigr) \, dx
  \end{split}
\end{equation*}

Consequently, \eqref{c1.2} yields 
\begin{equation}
2 \sum_{k=3}^\infty \frac{C(k,k)}{k^{k+1}}= \z(3)-3\z(2)-
\int_{x=0}^\infty \log\bigl(1-(1+x)e^{-x}\bigr) \, dx.
\end{equation}

Together with \eqref{A5} and \eqref{xx1}, this verifies (b).

{\bf Case (c):}  $1\leq k\leq n$, $j\ge1$ (Complex components).

Let 
\begin{equation}
  \kkc(p)=\kkcn(p)\=\sum_{k=1}^\infty\sum_{j=1}^\infty\kk(k,j,p),
\end{equation}
\ie, the number of complex components in $\gnp$, and
\begin{equation}\label{fn}
  f_n(p)=\E\kkc(p)=\sumk\sumj\E\k(k,j,p),
\end{equation}
the expected number of complex components in $G_{n,p}$.
The contribution to \eqref{E2} from the complex components is thus 
$\intp f_n(p)\,dp$. We make a change of variables and let
\begin{equation}\label{pgl}
  p=n^{-1}+\gl n^{-4/3},
\end{equation}
which means that we focus on the critical window. We will assume this
relation between $p$ and $\gl$ in the rest of the proof. We thus define
$\ff_n(\gl)=f_n(p)=f_n(\pn)$, and obtain the contribution,
letting $\ett{\dots}$ denote the indicator of an event,
\begin{equation}\label{em}
  \begin{split}
  \intp f_n(p)\,dp
&= 1-\frac1n + \intp\bigpar{f_n(p)-\ett{p>1/n}}\,dp
\\&
= 1-\frac1n + n^{-4/3}\intgl\bigpar{\ff_n(\gl)-\ett{\gl>0}}\,d\gl .	
  \end{split}
\end{equation}
We begin by showing that the integrand in the final integral converges
pointwise. We define, \cf{} \eqref{mgf2},
\begin{equation}\label{psi2}
  \psi_2(t)=\sum_{\ell=2}^\infty \wl t^l =\psi(t)-1-\sqrt{\xfrac\pi8}\,t,
\end{equation}
and
\begin{equation}\label{F}
  F(x,\gl)=\frac16x^3-\frac12x^2\gl+\frac12x\gl^2
=\frac{x}2\Bigpar{\gl-\frac x2}^2+\frac1{24}x^3.
\end{equation}

\begin{sublemma}\label{SLA}
  For any fixed $\gl\in\oooo$, as \ntoo,
  \begin{equation}\label{sla}
\ff_n(\gl)\to f(\gl)
=\frac1{\sqrt{2\pi}}\intx \psi_2(x^{3/2})e^{-F(x,\gl)}x^{-5/2}\,dx.	
  \end{equation}
\end{sublemma}

\proofstart
We note first that the integral in \eqref{sla} is convergent; 
for small $x$ we have $\psi_2(x)=O(x^2)$ and for large $x$  we have
$\psi_2(x)=O \bigpar{x^2e^{x^2/24}}$ by \eqref{jesp} while 
$e^{-F(x,\gl)} \le e^{-x^3/6+\gl x^2/2}=O(e^{-x^3/7})$ by \eqref{F}, remember that $\l$ is fixed in the integral.

We convert the sum over $k$ in \eqref{fn} 
to an integral by setting $k=\kx$.
Thus
\begin{equation}\label{fn2}
\ff_n(\gl)=
  f_n(p)=\intx\sumj\E\k\bigpar{\kx,j,p}n\qqqb\,dx.
\end{equation}
For any fixed $\l$ and fixed $x>0,\,j\ge1$, 
and $p=\pn$ and $k=\kx$ as above,
we have as \ntoo{}
by \eqref{ekk} and \eqref{wright} and standard calculations, see \eg{}
\cite[Section 4]{SJ179} 
or \cite[Section 11.10]{AS}
for further details,
\begin{equation*}
  \begin{split}
  \E\k(k,j,p)
&\sim
\frac{n^k}{k!}\exp\Bigpar{-\frac{k^2}{2n}-\frac{k^3}{6n^2}}
C(k,k+j)n^{-k-j}\bigpar{1+\gl n^{-1/3}}^k\exp\bigpar{-p(nk-k^2/2)}
\\&
\sim	
n^{-j}\frac{C(k,k+j)}{k!}\exp\bigpar{-k-F(k n^{-2/3},\gl)}
\\&
\sim (2\pi)\qqw \wj	
k\qw\Bigparfrac{k^{3/2}}{n}^j
e^{-F(k n^{-2/3},\gl)}
\\&
\sim n\qqqbw(2\pi)\qqw \wj	
 x^{3j/2-1}
e^{-F(x,\gl)}.
  \end{split}
\end{equation*}
Thus, as \ntoo,
\begin{equation}
\label{mag}
 n\qqqb \E\k(\kx,j,p)\to (2\pi)\qqw\wj x^{3j/2-1}e^{-F(x,\gl)}.
\end{equation}

Moreover, Bollob\'as \cite[Theorem 5.20]{B} has shown the uniform bound
\begin{equation}\label{ukk}
  C(k,k+j) \le \parfrac{\CC}{j}^{j/2} k^{k+(3j-1)/2} \CCdef\CCB
\end{equation}
for some constant $\CCB$ and all $k,j\ge1$.
Let $A\ge1$ be a constant, and consider first only components of size
$k\le A n^{2/3}$.
For such $k$, all $j\ge1 $ and $p=\pn$,
\eqref{ekk} and \eqref{ukk} yield by calculations similar to 
those above,
\begin{equation*}
  \begin{split}
  \E\k(k,j,p)
&\le \CC
\frac{n^k}{k!}\exp\Bigpar{-\frac{k^2}{2n}}
C(k,k+j)n^{-k-j}\bigpar{1+\gl n^{-1/3}}^{k+j}\exp\bigpar{-p(nk-k^2/2-j)}
\\
&\le \CC
n^{-j}\frac{C(k,k+j)}{k!}e^{-k+j\times o(1)}
\\
&\le\CCx n^{-j}
\Bigparfrac{2\CCB}{j}^{j/2}
{k^{3j/2-1}}
  \end{split}
\end{equation*}
(with $C_3$ possibly depending on $A$) and thus
\begin{equation*}
n\qqqb  \E\k(k,j,p)
\le\CCx 
\parfrac{\CC A^{3/2}}{j}^{j/2}.
\end{equation*}
The sum over $j$ of the right-hand side converges, and thus 
\eqref{mag} and dominated convergence yield, recalling \eqref{psi2},
\begin{equation}\label{sw}
\int_{x=0}^A\sumj\E\k\bigpar{\kx,j,p}n\qqqb\,dx
\to 
\frac1{\sqrt{2\pi}}\int_{x=0}^A \psi_2(x^{3/2})e^{-F(x,\gl)}x^{-5/2}\,dx.	
\end{equation}

For $k> An\qqqb$ we use the fact shown in \cite[(6.6)]{SJ179} that the
expected number of vertices in tree components of size at most $n\qqqb$ is
$n-O(n\qqqb)$; consequently, the expected number of vertices in all
components (complex or not) of size larger than $n\qqqb$ is $O(n\qqqb)$, and
the expected number of components larger than $An\qqqb$ is $\le \CC/A$.
The left-hand side of \eqref{sw} thus converges uniformly to the right-hand
side of \eqref{fn2} as \ntoo, and the result \eqref{sla} follows from
\eqref{sw} by letting $A\to\infty$.
\proofend

The next step is to use dominated convergence in \eqref{em}. For this we use
the following estimates.
For convenience, we let 
$\kkc(\pn)$ and its expectation
$\ff_n(\gl)$ 
be defined for all
real $\gl$, by  trivially defining $\kkc(p)=\kkc(0)=0$ for $p<0$ and
$\kkc(p)=\kkc(1)=1$ for $p>1$.

\begin{sublemma}\label{SLB}
  There exist integrable functions $g_1(\gl),g_2(\gl),g_3(\gl)$, 
not depending on $n$,  
such that
\begin{romenumerate}
\item 
  \begin{equation*}
	\ff_n(\gl)=\E \kkc(\pn) \le g_1(\gl),
\qquad \gl\le0,
  \end{equation*}
\item 
  \begin{equation*}
	\P\bigpar{\kkc(\pn)=0} \le g_2(\gl),
\qquad \gl\ge0,
  \end{equation*}
\item 
  \begin{equation*}
\ff_n(\gl)-1=
	\E \kkc(\pn)-1 \le g_3(\gl),
\qquad \gl\ge0.
  \end{equation*}
\end{romenumerate}
\end{sublemma}

\proofstart
  We use the method in Janson \cite{SJmulticyclic}.
We consider $G(n,p)$,
  $p\in\oi$,  
as a random graph process in the usual way: 
we  regard $p$ as time, edges are
added as $p$ grows from 0 to 1, and an edge $e$ is added at a time $T_e$ 
with a uniform distribution on $\oi$, with all $T_e$ independent.

As $G(n,p)$ evolves, there are at first only tree components, but later
unicyclic components and complex components appear as edges are added to the
graph. 
If we consider the complex components only, a new complex component is
created if a new edge is added to a unicyclic component, or if it joins two
unicyclic components. (Note that these are the only possibilities;
we do not regard the growth of an already existing complex component
as creating a new complex component.
Creation of a new complex component may happen one or several times. 
It is 
shown in \cite{SJbirth} that it happens only once with probability
converging to $5\pi/18$, but we will not need this.)
As evolution continues, the complex components may grow by merging with
trees or unicyclic components, and they may merge with each other, until at
the end only one complex component remains, containing all vertices. 

Let $\gf_n(k,p)$ be the intensity of creation on new complex components of
size $k$, \ie, the probability of creating a new complex component of size
$k$ in the interval $[p,p+dp]$ is $\gf_n(k,p)\,dp$. (For $p<0$, $p>1$ or
$k>n$, we  set $\gf_n(k,p)=0$.)
Further, let $\Phi_n(p)=\sum_{k\ge1}\gf_n(k,p)$, the intensity of creation
of complex components regardless of size.
We change variables as above and define also
\begin{align*}
  \psi_n(x,\gl)&=n\qqqbw \gf_n(\kx,\pn),
\\
\Psi_n(\gl)&=n^{-4/3} \Phi_n(\pn) = \intx \psi_n(x,\gl)\,dx.
\end{align*}
(The notation is not exactly as in \cite{SJmulticyclic}, where
the two ways of creating a complex component are treated separately, but the
estimates are the same.)

We have 
\begin{equation*}
  \gf_n(k,p)=\binom nk \hC(k) p^k(1-p)^{(n-k)k+\binom k2 -k-1}
\end{equation*}
where $\hC(k)$ is the number of ways to create a multicyclic component by
either adding an edge to a unicyclic component on $[k]$ or adding an edge
joining two unicyclic components whose vertex sets are complementary subsets
of $[k]$. The first case contributes
\begin{equation*}
  C(k,k)\lrpar{\binom k2 -k}=O(k^{k+3/2})
\end{equation*}
to $\hC(k)$ and the second
\begin{equation*}
\frac12\sum_{i=3}^{k-3}\binom ki C(i,i) C(k-i,k-i) i(k-i)
\le \CC \sum_{i=3}^{k-3} \binom ki e^i i!\, e^{k-i} (k-i)!
\le\CCx k e^k k!
=O(k^{k+3/2});
\end{equation*}
hence
\begin{equation*}
  \hC(k)=O(k^{k+3/2})=O\bigpar{ k e^k k!}.
\end{equation*}
(Cf.\ the more precise \cite[(2.30)]{SJmulticyclic}.)
The intensity $\psi_n(x,\gl)$ is bounded in 
\cite[(2.12)--(2.19)]{SJmulticyclic} by calculations similar to those in the
proof of Sublemma \ref{SLA}.
(In these bounds, and our versions below, 
$\gd,\gd_1,\dots$ are some positive constants.)

We use the results of \cite{SJmulticyclic} with some small modifications:
Equation (2.12) of \cite{SJmulticyclic} 
shows (together with the comments after it) 
that 
$$\psi_n(x,\gl)\leq \CC x e^{-\gd x^3-\gd x\gl^2}\qquad\text{ for }k\leq\d_1n
\text{ and }
-n^{1/3}\leq \l\leq \d_2n^{1/3}.
\CCdef\CCmci
$$ 
Then one line before (2.15) of \cite{SJmulticyclic} proves that
$$
\psi_n(x,\gl)\leq \CC xe^{-\d x^3-\d_3x\l n^{1/3}/3}\qquad\text{ for }k\leq
\d_3n\text{ and }\l\geq \d_2n^{1/3}.
\CCdef\CCmcii
$$
Because $\l\leq n^{4/3}$ always, 
it is legitimate to replace $-\d_3x\l n^{1/3}/3$ by
$-\d_3x\l^{5/4}$ to give 
$$\psi_n(x,\gl)\leq \CCmcii xe^{-\d x^3-\d_3x\l^{5/4}/3}\qquad\text{ for }k\leq \d_3n\text{ and }\l\geq \d_2n^{1/3}.$$
Then (2.17) of \cite{SJmulticyclic} proves that
$$\psi_n(x,\gl)\leq \CC ne^{-2\d_5n}
\qquad\text{ for }
\min\set{\d_1,\d_3}n\leq k\leq n.$$
We replace this by, using $\min\set{\d_1,\d_3}n^{1/3}\leq x\leq n^{1/3}$ and
$\l\le n^{4/3}$, 
$$\psi_n(x,\gl)\leq \CC x e^{-\gd_5 x^3}(1+\gl^4)^{-1}.
\CCdef\CCmciii
$$
We therefore have, for all $x$ and $\l$ 
(recalling that $\psi_n(x,\l)=0$ if $x>n^{1/3}$, $\gl<-n^{1/3}$ or 
$\l>n^{4/3}$), 
\begin{equation}\label{g}
  0 \le \psi_n(x,\gl) \le g(x,\gl)
=
\CCmci x e^{-\gd x^3-\gd x\gl^2}
+ \CCmcii x e^{-\gd x^3- \d_3x|\gl|^{5/4}/3}
+ \CCmciii x e^{-\gd_5 x^3}(1+\gl^4)\qw.
\end{equation}

Integrating we find
\begin{equation}\label{Psibound}
  \Psi(\gl) \le \intx g(x,\gl)\,dx 
\le \frac{\CC}{1+|\gl|^{5/2}}. \CCdef\CCPsi
\end{equation}

The number of complex components at any time is at most the number of
complex components that have been created so far. Taking expectations we
thus obtain, using \eqref{Psibound},
\begin{equation}
  \ff_n(\gl)=\E\kkc(\pn) \le \int_{\mu=-\infty}^\gl\Psi(\mu)\,d\mu
\le \int_{\mu=-\infty}^\gl \frac{\CCPsi}{1+|\mu|^{5/2}}\,d\mu .
\end{equation}
This verifies (i), with $g_1(\gl)=\CC(1+|\gl|^{3/2})\qw$ for $\gl\le0$.

Similarly, if there is no complex component at some time, at least one
complex component has to be created later. Thus,
\begin{equation}
\P\bigpar{\kkc(\pn)=0} 
\le \int_{\mu=\gl}^{\infty} \Psi(\mu)\,d\mu
\le \int_{\mu=\gl}^\infty \frac{\CCPsi}{1+|\mu|^{5/2}}\,d\mu ,
\end{equation}
which verifies (ii) with 
$g_2(\gl)=\CC(1+\gl^{3/2})\qw$ for $\gl\ge0$.

For (iii), let 
$Y(p)=\binom{\kkc(p)}2$ be the number of pairs of  complex components in
$\gnp$. Since $\kkc(p)-1\le Y(p)$, it suffices to estimate $\E Y(p)$.

If there is a pair of complex components in $\gnp$, then these components
have been created at some times $p_1$ and $p_2$ with $p_1\le p_2\le p$. The
intensity of this happening, with sizes $k_1=\ceil{x_1n\qqqb}$
and $k_2=\ceil{x_2n\qqqb}$ of the components at the moments of their
creations, is bounded in 
\cite[(2.24)--(2.26)]{SJmulticyclic} by (using modifications as above, and 
$g$ is defined in \eqref{g}),
\begin{equation*}
  \CC g(x_1,\gl_1)g(x_2,\gl_2)\,d\gl_1\,d\gl_2\,dx_1\,dx_2. \CCdef\CCsw
\end{equation*}
Moreover, if the two components still are distinct components in $\gnp$,
then, at least (ignoring further conditions from the growth of the
components), the original vertex 
sets of sizes $k_1$ and $k_2$ are not connected by any edge in the time
interval $[p_2,p]$; the (conditional) probability of this is
\begin{equation*}
\Bigpar{ 1-\frac{p-p_2}{1-p_2}}^{k_1k_2}
\le
\xpar{ 1-(p-p_2)}^{k_1k_2}
\le e^{-k_1k_2(p-p_2)}
\le e^{-x_1x_2(\gl-\gl_2)}.
\end{equation*}
Consequently, 
\begin{equation*}
  \begin{split}
  \ff_n(\gl)-1 
&\le \E Y(\pn)
\\&
\le g_3(\gl)
= \int_{\gl_1=-\infty}^{\gl} \int_{\gl_2=\gl_1}^{\gl} 
\int_{x_1=0}^\infty \int_{x_2=0}^\infty
\CCsw g(x_1,\gl_1)g(x_2,\gl_2)e^{-x_1x_2(\gl-\gl_2)}\,d\gl_1\,d\gl_2\,dx_1\,dx_2. 	
  \end{split}
\end{equation*}
This yields (iii), but it remains to verify that 
$\int_{\gl=0}^\infty g_3(\gl)\,d\gl<\infty$. 
Indeed, by Fubini and \eqref{g},
\begin{equation*}
  \begin{split}
&\intgloo g_3(\gl)\,d\gl
\\&\quad
= \int_{\gl_1=-\infty}^{\infty} \int_{\gl_2=\gl_1}^{\infty} 
\int_{x_1=0}^\infty \int_{x_2=0}^\infty
\CCsw g(x_1,\gl_1)g(x_2,\gl_2)
\int_{\gl=\gl_2}^\infty 
 e^{-x_1x_2(\gl-\gl_2)}\,d\gl\,d\gl_1\,d\gl_2\,dx_1\,dx_2 	
\\&\quad
= \int_{\gl_1=-\infty}^{\infty} \int_{\gl_2=\gl_1}^{\infty} 
\int_{x_1=0}^\infty \int_{x_2=0}^\infty
\CCsw \frac{g(x_1,\gl_1)g(x_2,\gl_2)}{x_1x_2}
\,d\gl_1\,d\gl_2\,dx_1\,dx_2 	
\\&\quad
\le\CCsw\lrpar{\int_{\gl=-\infty}^{\infty}\int_{x=0}^\infty  
\frac{g(x,\gl)}{x}
\,d\gl\,dx}^2
<\infty.
  \end{split}
\end{equation*}
\proofend

Sublemma \ref{SLB}(ii) implies that $1-\ff_n(\gl)\le g_2(\gl)$ for $\gl\ge0$,
and thus Sublemma \ref{SLB} yields
\begin{equation*}
  \bigabs{\ff_n(\gl)-\ett{\gl>0}} \le
	\begin{cases}
	  g_1(\gl), & \gl\le 0, \\
g_2(\gl)+g_3(\gl), & \gl>0.
	\end{cases}
\end{equation*}
This justifies using dominated convergence in the integral in \eqref{em},
and
Sublemma \ref{SLA} implies
\begin{equation}\label{per}
  \intgl\bigpar{\ff_n(\gl)-\ett{\gl>0}}\,d\gl 
\to
\ccc=  \intgloo\bigpar{f(\gl)-\ett{\gl>0}}\,d\gl .
\end{equation}
Hence \eqref{em} yields
\begin{equation}\label{em2}
  \begin{split}
  \intp f_n(p)\,dp
= 1-\frac1n + \ccc n^{-4/3}+o(n^{-4/3}) ,	
  \end{split}
\end{equation}
which is the sought result except for the expression for $\ccc$.

We transform the expression for $\ccc$ in \eqref{per} by first writing it as
\begin{equation}\label{er}
  \begin{split}
  \ccc&=\lim_{A\to\infty} \lrpar{-A+\int_{\gl=-\infty}^A f(\gl)\,d\gl}
\\&
=\lim_{A\to\infty} \lrpar{-A+ \frac1{\sqrt{2\pi}}\int_{\gl=-\infty}^A 
\intx \psi_2(x^{3/2})e^{-F(x,\gl)}x^{-5/2}\,dx\,d\gl}.	
  \end{split}
\end{equation}
By \eqref{psi2} we have $\psi_2(t)=O(t^2)$ for small $t$, which together
with \eqref{jk} shows that
\begin{equation*}
  \psi_2(t) = O\bigpar{t^2 e^{t^2/24}},
\qquad t\ge0
\end{equation*}
and thus by \eqref{F}, for all $x>0$ and $\gl\in\oooo$,
\begin{equation*}
  \psi_2(x^{3/2}) e^{-F(x,\gl)}
\le \CC x^3 e^{-x(\gl-x/2)^2/2}.
\end{equation*}
Hence, for $A>0$, with the substitutions $x=2A+s$ and $\gl=A-t$,
\begin{equation*}
  \begin{split}
\int_{x>2A}\int_{\gl< A} & \psi_2(x^{3/2})e^{-F(x,\gl)}x^{-5/2}\,dx\,d\gl
\le \CCx 
\int_{x>2A}\int_{\gl< A} e^{-x(\gl-x/2)^2/2}x^{1/2}\,dx\,d\gl
\\&
=
\CCx 
\int_{s>0}\int_{t>0} e^{-(2A+s)(t+s/2)^2/2}(2A+s)^{1/2}\,dt\,ds
\\&
\le
\CCx 
\int_{s>0}\int_{t>0} e^{-(2A+s)(t^2/2+s^2/8)}(2A+s)^{1/2}\,dt\,ds
\\&
=
\CC
\int_{s>0} e^{-(2A+s)s^2/8}\,ds
\le \CC A\qqw.
  \end{split}
\end{equation*}
Similar estimates show also
\begin{equation*}
  \begin{split}
\int_{x<2A}\int_{\gl> A} & \psi_2(x^{3/2})e^{-F(x,\gl)}x^{-5/2}\,dx\,d\gl
\le \CC \int_{s=0}^{2A} e^{-(2A-s)s^2/8}\,ds
\le \CC A\qqw.
  \end{split}
\end{equation*}
Consequently, we can subtract and add these integrals to \eqref{er},
yielding
\begin{equation}\label{er2}
  \ccc
=\lim_{A\to\infty} \lrpar{-A+ \frac1{\sqrt{2\pi}}\int_{\gl=-\infty}^\infty 
\int_{x=0}^{2A} \psi_2(x^{3/2})e^{-F(x,\gl)}x^{-5/2}\,dx\,d\gl}.
\end{equation}
It follows from \eqref{F} that
\begin{equation}
  \intgloo e^{-F(x,\gl)}\,d\gl
=e^{-x^3/24}\intgloo e^{-x(\gl-x/2)^2/2}\,d\gl
=e^{-x^3/24} \sqrt{2\pi/x}.
\end{equation}
Hence \eqref{er2} yields by Fubini
\begin{equation}\label{er3}
  \ccc
=\lim_{A\to\infty} \lrpar{-A+ 
\int_{x=0}^{2A} \psi_2(x^{3/2})e^{-x^3/24}x^{-3}\,dx}
=
\int_{x=0}^{\infty} \Bigpar{x^{-3}\psi_2(x^{3/2})e^{-x^3/24}-\frac12}\,dx.
\end{equation}

This completes the proof of Lemma \ref{lem1} and the proof of Theorem \ref{th1}.
\proofend

\section{Final remarks}\label{Sfinal}

\begin{remark}
We have shown that when the $X_e$ are uniform $[0,1]$ then $\E(L_n)$
converges to $\z(3)$ with an error term of 
order $1/n$. The constant $c_1$ is positive
and so for
large $n$ we have $\E(L_n)>\z(3)$. 
Fill and Steele \cite{Ste2} computed $\E(L_n)$ for $n\leq 8$.
$\E(L_n)$ increased monotonically and  it was 
natural to conjecture from this that $\E(L_n)$ increases monotonically for
all $n$. However, since $\E(L_n)$ converges to $\zeta(3)$ from above, we now
see that this turns out not to be true.
Note, however, that $c_2<0$, and that $|c_2|$ is much larger than
$c_1$. Thus we expect that $\E L_n>\zeta(3)$ only for very large $n$.

We have, if our numerical estimates are correct, $|c_2|/c_1\approx 45$, so
  a naive guess, ignoring higher order terms, would be that $\E L_n>\zeta(3)$
for $n > 45^3 \approx 10^5$.
We don't want to conjecture this, as we have no idea about the next term.
\end{remark}

\begin{remark}  
  By \eqref{eqLZ}, we obtain for $\E_{exp}(L_n)$ the same result as in
  \refT{th1} except that $c_1$ is increased by $\zeta(3)$ (while $c_2$
  remains the same). This gives a somewhat simpler $c_1$, which suggests
  that this version might be slightly simpler to analyze. Note that the
  formula \eqref{E2} holds for $\E_{exp}(L_n)$ if we replace $G_{n,p}$
  by the multigraph where each pair of vertices is connected by a $\Po(t)$
  number of edges, and integrate for $t\in(0,\infty)$. This suggests that it
  might be profitable to make a version of the argument below using these
  multigraphs, but we have not pursued this. (Cf.\ the use of multigraphs in
  \cite{SJbirth}.)   
\end{remark}

{\bf Acknowledgement:} In an earlier version of this paper, we showed that  
$\E(L_n)=\z(3)+\frac{c_1+o(1)}{n}$. Nick Read pointed us to his article \cite{read} 
which suggested that the $o(1/n)$ term could be replaced
by $(c+o(1))/n^{4/3}$. This encouraged us to go the extra mile and find the next term
and prove Nick's conjecture.

\appendix
\section{Appendix: Estimation of $c_2$}
The constant $c_2$ in \refT{th1} is the sum of the three coefficients for
$n^{-4/3}$ in \refL{lem1}(a)--(c), which we denote by $\cca$, $\ccb$ and
$\ccc$. 
By the change of variable $t=x^3/24$, 
and integration by parts (\cf{} \cite[\S5.9.5]{NIST}),
we obtain, as said in \refR{Rc2ab}, 
\begin{align}
\cca&=
\frac{24^{-2/3}}{3}\int_{t=0}^\infty t^{-5/3}(e^{-t}-1)\,dt
=-\frac18 3^{-2/3} \Gamma\Bigpar{\frac13}
= - 0.16098\dots, \label{cca}
\\
\ccb&=
\sqrt{\frac{\pi}8}\frac{24^{-1/6}}{3}\int_{t=0}^\infty t^{-7/6}(e^{-t}-1)\,dt
=-\frac123^{-1/6}\sqrt\pi\,\Gamma\Bigpar{\frac56} = - 0.83298\dots. \label{ccb}
\end{align}

The coefficient $\ccc$ is given by an integral in \refL{lem1}, see also 
\eqref{er3}.
To evaluate $\ccc$, we  change variables by $x=y\qqq$
and use the definition \eqref{psi2}
of $\psi_2$ to obtain
\begin{equation}
  \begin{split}
\ccc&=
\frac13
\int_{y=0}^{\infty} \Bigpar{y^{-1}\psi_2(y^{1/2})-\frac12e^{y/24}}e^{-y/24}
 y^{-2/3}\,dy
\\&
=
\frac13
\int_{y=0}^{\infty} 
\sum_{k=1}^\infty
\Bigpar{w_{2k}y^{k-1}+w_{2k+1}y^{k-1/2}-\frac{y^{k-1}}{2\cdot24^{k-1}(k-1)!}}
e^{-y/24} y^{-2/3}\,dy.
  \end{split}
\end{equation}
We interchange the order of integration and summation, which is justified
below, and obtain
\begin{equation}\label{pil}
  \begin{split}
\ccc&=
\frac13
\sum_{k=1}^\infty
\int_{y=0}^{\infty} 
\Bigpar{w_{2k}y^{k-1}+w_{2k+1}y^{k-1/2}-\frac{y^{k-1}}{2\cdot24^{k-1}(k-1)!}}
e^{-y/24} y^{-2/3}\,dy
\\&
=
\frac{24\qqq}{3}
\sum_{k=1}^\infty
\Bigpar{w_{2k}24^{k-1}\Gamma(k-2/3)+w_{2k+1}24^{k-1/2}\Gamma(k-1/6)
 -\frac{\Gamma(k-2/3)}{2\,\Gamma(k)}}.
  \end{split}
\end{equation}
We note that \eqref{wbex} and \eqref{ebex} yield, together with Stirling's
formula, 
$\wl\sim6\cdot 24^{-\ell/2}/\Gamma(\ell/2)$, which implies that
\begin{equation*}
w_{2k}24^{k-1}\Gamma(k-2/3)\sim w_{2k+1}24^{k-1/2}\Gamma(k-1/6)
\sim \tfrac14 k^{-2/3}
\qquad \text{as $k\to\infty$}
\end{equation*}
so the three terms in the sum in \eqref{pil} are all of order $k^{-2/3}$,
showing that we cannot sum them separately.
However, their leading terms cancel. 
A more precise calculation using \eqref{fbex} yields
\begin{equation}
  \E \bex^r = \sqrt{18}\,r\Bigparfrac{r}{12e}^{r/2}\bigpar{1+O(r\qw)},
\qquad r>0,
\end{equation}
and thus by \eqref{wbex} and Stirling's formula,
\begin{equation}
 \wl = \frac{3\sqrt\ell}{\sqrt{\pi}}\Bigparfrac{e}{12\ell}^{\ell/2}
\bigpar{1+O(\ell\qw)}
= \frac{6\cdot 24^{-\ell/2}}{\Gamma(\ell/2)}
\bigpar{1+O(\ell\qw)},
\qquad \ell\ge1.
\end{equation}
Hence, 
\begin{equation}\label{sos}
w_{2k}24^{k-1}\Gamma(k-2/3)
= \tfrac14 k^{-2/3} \bigpar{1+O(k\qw)},
\qquad \text{as $k\to\infty$},
\end{equation}
and the same estimate holds for $w_{2k+1}24^{k-1/2}\Gamma(k-1/6)$,
while $\Gamma(k-2/3)/\Gamma(k)=k^{-2/3} \bigpar{1+O(k\qw)}$.
Consequently, the summand in \eqref{pil} is $O(k^{-5/3})$.

The constants $w_k$ can be computed by a recursion formula, see
\cite{Wright} and \cite{SJ201}, and a numerical summation of the first 1000
terms in \eqref{pil} yields $-0.7331$.  
It can be shown, 
using \eqref{fbex} with the further second order term 
given in \cite{SJ203}
(which replaces $O(x^{-2})$ by $-\frac19x^{-2}+O(x^{-4})$),
that the terms in the sum in \eqref{pil} are
$\sim -\frac16 k^{-5/3}$, and using this to
estimate the sum of the terms with $k>1000$
yields the estimate
$\ccc\approx -0.7355$   
which together with \eqref{cca}--\eqref{ccb} yields
\begin{equation}\label{A4}
  c_2\approx -1.7295.   
\end{equation}
{\em The tail estimate is not rigorous. Replacing $O(x^{-4})$ by $\leq Cx^{-4}$ for some estimate $C$ is what is needed
to make the tail estimate rigorous. Nevertheless, it seems unlikely that the estimate in \eqref{A4} is very far off.}

To justify the interchange of summation and integration above, it is by
Fubini's theorem sufficient to verify that
\begin{equation}\label{pilk}
\sum_{k=1}^\infty
\int_{y=0}^{\infty} 
\lrabs{w_{2k}y^{k-1}+w_{2k+1}y^{k-1/2}-\frac{y^{k-1}}{2\cdot24^{k-1}(k-1)!}}
e^{-y/24} y^{-2/3}\,dy
<\infty.
\end{equation}
Indeed, we claim that the integral in \eqref{pilk} is $O(k^{-7/6})$.
Using \eqref{sos}, its analogue for $2k+1$, and
$\Gamma(k-2/3)/\Gamma(k)=k^{-2/3}(1+O(k\qw))$, it follows easily that the
integral is, after another change of variable $t=y/24$, 
\begin{equation}\label{ull}
\frac{24\qqq}{4}k\qqqbw
\int_{t=0}^{\infty} 
\lrabs{\frac{t^{k-7/6}}{\Gamma(k-1/6)}-\frac{t^{k-5/3}}{\Gamma(k-2/3)}}
e^{-t}\,dt + O\bigpar{k^{-5/3}}.
\end{equation}
Let $I_k$ denote the integral in \eqref{ull}. By the Cauchy--Schwarz
inequality,
\begin{equation*}
  \begin{split}
I_k^2 
&\le
\int_{t=0}^{\infty} t^{k-1} e^{-t}\,dt \cdot
\int_{t=0}^{\infty} 
\biggpar{\frac{t^{k-7/6}}{\Gamma(k-1/6)}-\frac{t^{k-5/3}}{\Gamma(k-2/3)}}^2
t^{1-k}e^{-t}\,dt 
\\&
= \gG(k)\lrpar{
\frac{\gG(k-2/6)}{\gG(k-1/6)^2}
-2 \frac{\gG(k-5/6)}{\gG(k-1/6)\gG(k-4/6)}
+ \frac{\gG(k-8/6)}{\gG(k-4/6)^2}}
\\&=O\bigpar{k\qw}.
	  \end{split}
\end{equation*}
Consequently, $I_k=O(k\qqw)$, which shows that \eqref{ull} is $O(k^{-7/6})$,
and thus \eqref{pilk} holds as claimed above.


\begin{thebibliography}{99}

\newcommand\vol{\textbf}
\newcommand\no[1]{\relax}

\bibitem{AS}
 N. Alon and J. H. Spencer,
 \emph{The Probabilistic Method}, 3rd ed,
 Wiley, New York (2008).


\bibitem{BFM}  A. Beveridge, A. M. Frieze and C. J. H. McDiarmid,
{\em Minimum length spanning trees in regular graphs}, Combinatorica 18
(1998) 311--333. 

\bibitem{B} B. Bollob\'as, {\em Random Graphs}, Second Edition, Cambridge University Press (2001).

\bibitem{Ste2} J. A. Fill and M. J. Steele, 
{\em Exact expectations of minimal spanning trees for graphs with random
edge weights},
Stein's Method and Applications,
Singapore Univ. Press, Singapore (2005)  169--180. 

\bibitem{Flax} A. Flaxman, {\em The lower tail of the random minimum spanning tree},
Electronic Journal of Combinatorics 14 (2007).

\bibitem{F1} A. M. Frieze, {\em On the value of a random minimum spanning tree problem},
{Discrete Applied Mathematics} 10 (1985) 47--56.

\bibitem{FM}  A. M. Frieze and C. J. H. McDiarmid, {\em On random minimum length
  spanning trees}, 
Combinatorica 9 (1989) 363--374.

\bibitem{FRT} A. M. Frieze, M. Ruszink\'o and L. Thoma, {\em A note on random
  minimum length spanning trees}, 
Electronic Journal of Combinatorics 7 (2000) R41.

\bibitem{Ga} D. Gamarnik, 
{\em The expected value of random minimal length spanning tree of a complete
graph}.
Proceedings of the Sixteenth Annual ACM-SIAM
Symposium on Discrete Algorithms (SODA 2005),  ACM, New York (2005) 700--704. 

\bibitem{SJmulticyclic}
S. Janson,
{\em Multicyclic components in a random graph process},
Random Structures and Algorithms 4 (1993, 71--84.

\bibitem{Ja}  S. Janson,
{\em The minimal spanning tree in a complete graph and a functional limit
theorem for trees in a random graph},
Random Structures and Algorithms 7 (1995) 337--355.

\bibitem{SJ201}
S. Janson,
{\em Brownian excursion area, Wright's constants in graph enumeration, 
and other Brownian areas},
Probability Surveys, 
3 (2007) 80--145.

\bibitem{SJISE} 
S. Janson and P. Chassaing, 
{\em The center of mass of the ISE and the Wiener index of trees}.
Electronic Comm. Probability 9 (2004) paper 20, 178--187.


\bibitem{SJbirth}
S. Janson, D.E. Knuth, T. {\L}uczak and B. Pittel,
{\em The birth of the giant component},
Random Structures and Algorithms 3 (1993) 233--358.

\bibitem{SJ203}
S. Janson and G. Louchard,
{\em Tail estimates for the Brownian excursion area and other Brownian areas},
Electronic Journal Probability 12 (2007) 1600--1632.


\bibitem{JLR} S. Janson, T. {\L}uczak and A. Ruci\'nski, {\em Random Graphs}, 
Wiley, New York (2000).

\bibitem{SJ179}
S. Janson and J. Spencer,
{\em A point process describing the component sizes in the critical window of the
random graph evolution}, 
Combinatorics Probability and Computing 16 (2007) 631--658. 

\bibitem{LZ} W. Li and X. Zhang, {\em On the difference of expected lengths of
  minimum spanning trees}, 
Combinatorics, Probability and Computing 18 (2009) 423--434.

\bibitem{Lou0}
G. Louchard,
{\em Kac's formula, L\'evy's local time and Brownian excursion}. 
Journal of Applied Probability 21 (1984) 479--499.

\bibitem{Lou}
G. Louchard,
{\em The Brownian excursion area: a numerical analysis}. 
Computational Mathematics and Applications 10 (1984) 413--417.
Erratum:
Computational Mathematics and Applications Part A 12 (1986) 375.

\bibitem{NIST}
\emph{NIST Digital Library of Mathematical Functions}.
\url{http://dlmf.nist.gov/}


\bibitem{Pe} M. Penrose, {\em Random minimum spanning tree and percolation on the
  $n$-cube}, Random Structures and Algorithms 12 (1998) 63--82.

\bibitem{read} N. Read, 
{\em Minimum spanning trees and random resistor networks in $d$ dimensions},
Physics Revue. E 72, 036114 (2005).

\bibitem{Re} A. R\'enyi, {\em Some remarks on the theory of trees}, 
Publ. Math. Inst. Hungar. Acad. Sci. 4 (1959) 73--85. 

\bibitem{Spencer}
J. Spencer, 
{\em Enumerating graphs and Brownian motion}.  
Communications in Pure and Applied Mathematics 50 (1997)  291--294. 

\bibitem{Ste1} M. J. Steele, {\em On Frieze's $\z(3)$ limit for lengths of
  minimal spanning trees}, Discrete Applied Mathematics 18 (1987) 99--103.

\bibitem{Wastlund} J. W\"astlund, {\em An easy proof of the $\zeta(2)$
limit in the random assignment problem}, Electronic Communications in
Probability 14 (2009) 261--269. 

\bibitem{Wright}
E. M. Wright,
{\em The number of connected sparsely edged graphs}. 
Journal of Graph Theory 1 (1977) 317--330.


\end{thebibliography}
\end{document}